\documentstyle[12pt]{article}
\input amssym.def
\input amssym.tex

\newcommand{\nc}{\newcommand}

\newcommand{\bla}{\phantom{bbbbb}}
\newcommand{\onebl}{\phantom{a} }

\newcommand{\beq}{\begin{equation}}
\newcommand{\eeq}{\end{equation}}
\newcommand{\beqst}{\begin{equation*}}
\newcommand{\eeqst}{\end{equation*}}
\newcommand{\barr}{\begin{array}}
\newcommand{\earr}{\end{array}}
\newcommand{\beqar}{\begin{eqnarray}}
\newcommand{\eeqar}{\end{eqnarray}}
\newtheorem{theorem}{Theorem}[section]
\newtheorem{corollary}[theorem]{Corollary}

\newtheorem{prop}[theorem]{Proposition}
\newtheorem{definition}[theorem]{Definition}
\newtheorem{remit}[theorem]{Remark}

\newcommand{\RR}{{\Bbb R }}
\newcommand{\CC}{{\Bbb C }}
\newcommand{\ZZ}{{\Bbb Z }}

\newcommand{\call}{{\mbox{$\cal L$}}}
\newcommand{\calm}{{\mbox{$\cal M$}}}

\newcommand{\calo}{{\mbox{$\cal O$}}}

\newcommand{\calr}{{\mbox{$\cal R$}}}

\def\a{\alpha}
\def\b{\beta}
\def\g{\gamma}

\def\z{\zeta}

\def\l{\lambda}

\def\G{\Gamma}

\def\L{\Lambda}

\nc{\xb}{{X_\beta} } 
\nc{\pr}{\partial}
\nc{\cbi}{C^{*,*} (G) }
\nc{\ctbi}{ {\tilde C}^{*,*}(G) } 
\nc{\YY}{Y}
\nc{\mc}{\alpha}
\nc{\bmc}{\bar{\alpha} } 
\nc{\lieg}{\frak{g} } 
\nc{\symf}{\omega} 
\nc{\inpr}[1]{{ \langle {#1} \rangle} }
\nc{\epsr}{{ \epsilon_R} } 
\nc{\zt}{{Z_t} }
\nc{\teta}{{ \tilde{\eta} } } 
\nc{\mg}{{ \calm}}
\nc{\Proof}{ \noindent{\em Proof:} } 

\nc{\diag}{\bigtriangleup}
\nc{\epm}{e}
\nc{\epmb}{{e_\beta}}
\nc{\proj}{{\rm pr} } 
\nc{\sigfo}{\sigma}
\nc{\sigf}{\sigma}
\nc{\tw}{ {\tilde{ \omega}} } 
\nc{\hop}{I}
\nc{\lrar}{\longrightarrow}

\newcommand{\renorm}{{ \setcounter{equation}{0} }}

\nc{\FF}{ {\Bbb F} } 

\begin{document}

\title{Symplectic Forms on Moduli Spaces of Flat Connections
on 2-Manifolds}

\author{Lisa C. Jeffrey \\
Mathematics Department \\
Princeton University \\
Princeton, NJ 08544, USA
\thanks{This material is based on work 
supported by the National Science Foundation under Grant No.
DMS-9306029.} }

\date{April 1994; revised April 1996} 

 \maketitle

\begin{abstract}
Let $G$ be a compact connected semisimple Lie group. We extend the techniques
of Weinstein [W] to give a construction in group cohomology 
of  symplectic forms $\omega$ on \lq twisted' moduli spaces of 
representations of the fundamental group $\pi$ of a 2-manifold
$\Sigma$ (the smooth analogues of ${\rm Hom} (\pi_1(\Sigma), G)/G$)
and on relative character varieties of fundamental 
groups of 2-manifolds.
We extend this construction to exhibit a symplectic form on the 
extended moduli space [J1]  (a Hamiltonian 
$G$-space  from which  these moduli spaces may
be obtained by symplectic reduction), and compute the moment map 
for the action of $G$ on the extended moduli space.
\end{abstract}

\renorm
\section{Introduction}

Let $\Sigma$ be a closed oriented 2-manifold of genus $g
\ge 2$; the fundamental 
group of $\Sigma$ will be denoted $\pi$. 
Let $G$ be a compact connected Lie group with Lie algebra
$\lieg$. This paper concerns the moduli space 
$\mg = {\rm Hom} (\pi, G)/G$ of conjugacy classes of representations
of $\pi$ into $G$, and certain more general analogues of $\mg$. 
The space $\mg$ has an open dense set on which the structure of 
a smooth symplectic manifold is defined.
In addition to the definition we have given in terms of representations
of $\pi$, the space $\mg$ has two alternative descriptions. The 
first of these is the gauge theory description: 
via the holonomy map, 
$\mg$ is identified with the space of gauge equivalence classes of flat 
connections on a trivial  principal 
$G$ bundle over $\Sigma$. The second alternative description appears
once one fixes a complex structure on the 2-manifold $\Sigma$, so that
$\Sigma$ becomes a Riemann surface; $\mg$ is then  identified with 
the space of equivalence classes of 
semistable holomorphic $G^\CC$ bundles over $\Sigma$.

The purpose of this paper is to extend work of Karshon [K] and 
Weinstein [W] to a more general setting; this paper follows [W]
closely and should be read in conjunction with it. These two papers 
complement the work of Goldman [G]. Goldman originally gave a construction
in group cohomology of the symplectic form $\symf$ on the space 
$\mg$; in order to prove $\symf$ is closed, Goldman used the gauge
theory description of the symplectic form. Karshon [K] gave a proof
of the closedness of the symplectic form using group cohomology; 
Weinstein [W] reinterpreted Karshon's construction in the setting 
of the de Rham-bar bicomplex [B,Sh]. In the present work we extend
Weinstein's work to construct symplectic forms 
on relative character varieties of 
surface groups, and
on `twisted' 
moduli spaces $\calm_\beta$ of bundles on Riemann surfaces (see Section 6 of 
[AB]) associated to an element $\beta \in Z(G)$. The spaces 
$\calm_\beta$ share many properties with $\calm$ but in 
general have less singularities. Indeed, $\calm_\beta$ is smooth when 
$G=SU(n)$ and $\beta$ is a generator of the center of 
$SU(n)$: in contrast, even when 
$G = SU(2)$, the space $\calm$ is smooth only in the very special case
$g = 2$.

We also give a group cohomology 
construction of symplectic forms on the 
extended moduli spaces $X$ and $\xb$  [J1]: these are 
 finite dimensional symplectic $G$-spaces
from which $\mg$ and $\calm_\beta$
 may be obtained by symplectic reduction.
 In [J1],
  symplectic structures on $\xb$  and $X$
(which are
in fact the same as the symplectic 
forms  we recover below)
 were specified using gauge theoretic techniques.
Here  we 
compute the moment maps for the $G$ action on $X$
and $\xb$; up to a normalization
factor, these coincide with the moment maps found in [J1]. 

The purpose of the construction 
of $\xb$ in [J1] was to exhibit $\calm_\beta$ as the result of 
finite dimensional symplectic reduction (in contrast to the infinite 
dimensional quotient construction given in [AB]). We  make further
use of this finite dimensional quotient construction 
in [JK2], where we  extend the techniques of [JK1] (which gives 
a formula for intersection pairings in the cohomology 
ring of the symplectic quotient $\calm_X$ of a finite dimensional 
Hamiltonian $G$-space $X$, in terms of the $G$-equivariant cohomology 
$H^*_G(X)$) to treat the intersection 
pairings in $\calm_\beta$, starting from the 
equivariant cohomology of $\xb$. By this means we  give proofs of 
formulas 
(found originally by Witten [Wi] using physical methods) for the 
intersection pairings in $H^*(\calm_\beta)$.

In this paper, 
the symplectic forms on $X$  and $\xb$ are
 constructed explicitly in 
terms of the Maurer-Cartan form on $G$ and   the chain 
homotopy operator that occurs in the standard proof of the 
Poincar\'e lemma. This explicit description of the symplectic form will be 
important in [JK2], where we make use of 
explicit
equivariantly closed differential forms representing the relevant classes
in de Rham cohomology.
In [J2] we extend the methods of this paper to give explicit 
representatives in De Rham cohomology for all the generators of the
cohomology ring of $\calm_\beta$ (one of which is the cohomology 
class of the symplectic form); for the  applications in [JK2], we 
shall need the de Rham representatives for all the generators.

After this work was completed, we 
received the paper of Huebschmann [H], in which he has obtained similar
 results independently.

{\em Acknowledgement:} We thank A. Weinstein for helpful conversations.

\renorm
\section{Group cohomology}

Let $\FF = \FF_{2g} $ be the free group on $2g$ generators
$x_1, \dots, x_{2g}$. We introduce a relation 
$R = \prod_{i = 1}^g [x_{2i-1}, x_{2i}], $ where $[a,b] $ 
denotes the commutator $ab a^{-1} b^{-1} $. The fundamental 
group $\pi$ of a closed 2-manifold of genus $g$ is then 
given by 
$$ \pi = \FF/\calr$$
where $\calr$ is the normal subgroup generated by $R$. 

We shall work with Eilenberg-Mac Lane group chains
(see [G] section 3.8), and shall denote the differential 
on the group chain complex by $\pr$. The $p$-chains 
$C_p(\G)$ on a group $\G$ are $\ZZ$-linear combinations 
of elements of $\G^p$. In particular, associated to the relation $R$ 
there is a distinguished 2-chain $c  \in C_2(\FF)$ given by 
\beq \label{1.1} c = \sum_{i =1}^{2g} (\partial R/\partial x_i, x_i). \eeq
(Here, $\partial/\partial x_i$ refers to the differential in the 
Fox free differential calculus: see [G], sections 3.1-3.3.) 
Goldman ([G], above  Proposition 3.9) shows 
\beq \label{1.2} \pr c = 1 - R. \eeq

\renorm
\section{The de Rham-bar bicomplex}

Weinstein [W] introduces a bicomplex $(\cbi,\delta, d)$ whose $p,q$ term
is $C^q(G^p)$. The second coboundary is the exterior 
differential $d$, while the first is the differential $\delta$ 
appearing in group  cohomology:
\beq (\delta \beta ) (g_0, \dots, g_p) = 
(-1)^{p+1} \beta (g_0, \dots,  g_{p-1} )
+ \sum_{i = 1}^p (-1)^i \beta (g_0, \dots, g_{i-1} g_{i}, \dots, g_p ) 
+  \beta(g_1, \dots, g_p). \eeq
Let $\YY$ denote ${\rm Hom} (\FF, G) = G^{2g}. $ Then there is a second
bicomplex $(\ctbi, \delta, d)$ whose $p,q$ term is 
$\Omega^q(\FF^p  \times Y) $ $ = \Omega^0 (\FF^p) \otimes \Omega^q(\YY)$. 
(The differential $\delta$ in the bicomplex $\ctbi$ is the adjoint 
of the differential $\partial $ in the Eilenberg-Mac Lane 
group chain complex.)
As in [W], the evaluation maps
$$E_p: \FF^p \times \YY \to G^p$$
give rise to maps 
$$E_p^*: \Omega^q(G^p) \to \Omega^0 (\FF^p) \otimes \Omega^q(\YY) $$
which combine to form a map of bicomplexes $E^*: \cbi \to \ctbi$. 

We recall the following elements of $\cbi$ defined in [W]. Let 
$\mc \in \Omega^1(G) \otimes \lieg$
denote the (left-invariant) Maurer-Cartan form on $G$, and 
$\bmc$ the corresponding right-invariant form. Define projection 
maps $\pi_i: G^p \to G$ ($i = 1, \dots, p$) to the $i$'th copy 
of $G$, and let $\mc_i = \pi_i^* \mc$ and 
$\bmc_i = \pi_i^* \bmc$. In terms of this notation, the following
are introduced in [W]:
\beq \l = \frac{1}{6} \mc \cdot [\mc, \mc] \onebl \in 
\Omega^3(G), \eeq
\beq \Omega = \mc_1 \cdot \bmc_2 \in \Omega^2 (G^2). \eeq
Given $\eta \in \lieg$, Weinstein also introduces 
\beq \theta_\eta = \eta \cdot (\mc + \bmc) \in \Omega^1(G). \eeq
(Here, $\cdot$ denotes an invariant inner product on $\lieg$.) 
These forms satisfy the following properties ([W], Lemmas
3.1, 3.3 , 4.1 and 4.4): 
\begin{prop} \label{p1} We have 
\beq \label{pe1} d\l = 0; \eeq
\beq \label{pe2} d \Omega = \delta \l; \eeq
\beq \label{pe4} \iota_\teta \l =  d \theta_\eta; \eeq
\beq \label{pe3} \iota_\teta \Omega = - \delta \theta_\eta, \eeq
where $\teta$ is the vector field  generated by $\eta$.

\end{prop}

\renorm
\section{Two-forms on moduli spaces}

Let us now introduce 
\beq   \label{4.1} \symf = \inpr{c, E^* \Omega} \in \Omega^2(\YY). \eeq
We then have 
\begin{prop} \label{p2} In terms of the identification of 
$\YY = {\rm Hom}(\FF, G) $ with $G^{2g}$, we have
\beq \label{dsymf} d \symf = - \epsr^* \l \eeq
where $\epsr: {\rm Hom}(\FF, G) \to G$  is the map  given by evaluation
on the
element $R$ $\in \FF$: 
\beq 
\label{epsrdef}
\epsr (g_1, \dots, g_{2g} ) = \prod_{i = 1}^g [g_{2i-1}, g_{2i}]. \eeq
\end{prop} 
\Proof We have
$$ d \inpr{c, E^* \Omega} = \inpr{c, d E^* \Omega} $$
$$ \bla = \inpr{c, E^* d \Omega} = \inpr{c, E^* \delta \l} $$ 
$$ \bla = \inpr{c, \delta E^*\l} = \inpr{\partial c, E^* \l} $$
$$ \bla = \inpr{1 - R, E^* \l}, $$  where the last 
step follows from (\ref{1.2}). $\square$

 If $t$ is any element of $G$, define $\YY_t = \epsr^{-1} (t) $ 
$\subset \YY$. The
following is an immediate consequence of Proposition \ref{p2}:
\begin{prop} \label{p3}
The 2-form $\symf$ restricts on $\YY_t$ to a closed form.
\end{prop}

\begin{definition} If $t$ is an element of $G$, let $\zt \subset G$ 
be the centralizer of $t$ in $G$. \end{definition}
Notice that $\YY_t$ carries an action of  $\zt$ by conjugation.

\begin{definition} 
If $t \in G$, the {\em relative character variety}
associated to $t$ is the space $\calm_t = \YY_t /\zt$, where
$\zt$ acts on $\YY_t$ by conjugation. \end{definition}
Some properties of the symplectic geometry of relative character
varieties were given in [JW1] and [JW2]. 
Relative character varieties also arise in algebraic geometry where  (under 
appropriate circumstances) they  are
identified with moduli spaces of semistable parabolic 
vector bundles on Riemann surfaces: see 
[MS] (where the identification between relative character varieties and
moduli spaces of parabolic bundles is given) or [Se].
A gauge theory construction of a symplectic form  on relative
character varieties is 
 given in [J1]: it will follow from the Remark at the end
of Section 5 that this is essentially the same as the symplectic form we 
construct here (i.e., the two are the same under a natural map 
identifying the relevant Zariski tangent spaces).

By extending the calculations in [W], we obtain the following two 
propositions. 

\begin{prop} \label{p4} The form $\symf$ is invariant under the
action of $G$ on $\YY$ by conjugation.
\end{prop}
\Proof Let $\eta \in \lieg$; we will show that  the Lie derivative
of $\symf$ 
with respect to the vector field $\teta$  generated by $\eta$ 
is zero, in other words 
$\call_\teta \symf = (d \iota_\teta + \iota_\teta d) \symf = 0 $.
Now 
\beq \iota_\teta d \symf =  - \iota_\teta \epsr^* \l = 
-\epsr^* \iota_\teta \lambda . \eeq
Also $d \iota_\teta \symf = d \inpr{c, E^* \iota_\teta \Omega}. $
But 
$\iota_\teta \Omega = - \delta \theta_\eta $ by (\ref{pe3}), so we get 
$$ d \iota_\teta \symf = - d \inpr{c, E^* \delta\theta_\eta } $$
$$ \bla = - d \inpr{c, \delta E^* \theta_\eta } $$ 
$$ \bla  = -  d \inpr{\partial c, E^* \theta_\eta} = -  d\inpr{1 - R, E^* 
\theta_\eta} $$ 
$$ =   d \epsr^* \theta_\eta = \epsr^* d \theta_\eta. $$ 
But by (\ref{pe4}) we have $ d \theta_\eta = \iota_\teta \l $, 
so $ d \iota_\teta \symf = \epsr^* \iota_\teta \l$, and 
 $\call_\teta \symf = 0 $. $\square$

\begin{prop} \label{p5} 
We have the following identification of 1-forms on 
$G^{2g}$: 
\beq \label{4.66}
\iota_\teta \symf  =   \epsr^* \theta_\eta. \eeq
Thus, the  restriction of $\symf$ to $\YY_t$ is horizontal
with respect to the action of $\zt$.
(In other words, if $\eta \in {\rm Lie} (\zt)$ generates the  vector
field $\tilde{\eta} $ on $\YY$, then $\iota_{\teta} \symf |_{\YY_t} = 0 $, 
where $\iota_{\teta} $ denotes the interior product with respect to 
$\teta$.)
\end{prop} 
\Proof We have $$\iota_\teta \symf = \inpr{c, \iota_\teta E^* \Omega} 
= \inpr{c, E^* \iota_\teta \Omega} = - \inpr{c, E^* \delta \theta_\eta},$$
$$ \onebl = - \inpr{\partial c, E^* \theta_\eta} = - \inpr{1 - R, 
E^* \theta_\eta} =   \epsr^* \theta_\eta. $$ 
This form necessarily restricts to zero on the level sets of $\epsr$.
$\square$

Propositions \ref{p4} and \ref{p5}  imply that the form 
$\symf$ descends to a 2-form $\bar{\omega}$ 
on the space $\calm_t$, which is 
closed by Proposition \ref{p3}. Nondegeneracy will be established 
in Corollary \ref{c5.5}.

In particular if $t $ is a central element  $\beta  
\in Z(G)$ then 
$Z_\beta$ is the full group $G$.  If $G = U(n)$ and 
$\beta$ 
$= e^{2 \pi i d/n} 
{\rm diag} (1 , \dots, 1)$ 
 then  the space $\calm_\beta$ 
appears in algebraic geometry  (see [AB]) as the moduli space
of semistable holomorphic vector bundles of rank $n$ and 
degree $d$. When $\beta$ is as above but 
 $G= SU(n)$, the algebraic geometry interpretation of 
the space 
$\calm_\beta$ is  as the moduli space
of semistable holomorphic vector bundles of rank $n$ and degree
$d$ with fixed determinant.\footnote{When  $G = U(n) $ or 
$SU(n)$ and 
$d$ is coprime to 
$n$, 
the spaces $\calm_\beta$ are
  smooth manifolds. This is in contrast to the spaces
$\mg$, which  are singular except in a few  
very special cases.}
 The constructions in this section  exhibit a group cohomology 
construction of a symplectic form on the \lq twisted' 
moduli spaces $\calm_\beta$
associated to central elements  $\beta$  of $G$.  It is shown in [J2]
that the form $\bar{\omega}$ is in the cohomology class of (a constant multiple
of) the standard generator $f_2$ of $H^*(\calm_\beta; \RR)$ 
(in the notation of Sections 2 and 9 of [AB]). In fact in [J2] 
we extend the construction which gives rise to the symplectic 
form $\bar{\omega}$,
 to give representatives in de Rham cohomology for all the 
generators of the ring $H^*(\calm_\beta; \RR)$ given in [AB].
Our applications in [JK2] will rely at least as heavily on the fact that
$\bar{\omega}$ is a de Rham representative of the cohomology class
$f_2$ as on its being nondegenerate:  in any case, many results of the 
type we shall invoke for 
Hamiltonian  $G$-manifolds generalize (cf. [KT]) to manifolds where 
the symplectic form degenerates on a locus of measure $0$.

\renorm
\section{The symplectic form on the extended moduli space}

Let $\beta$ be an element of the center $Z(G)$. 
The associated 
{\em extended moduli space}
 $X_\beta$ constructed in [J1] may be described
as a  fibre product 
\beq \xb = (\epsr \times \epmb)^{-1} (\diag) \subset G^{2g} \times \lieg. \eeq
Here, $\diag$ is the diagonal in $G \times G$ and 
$\epsr: G^{2g}  \to G$ was defined above, while $\epm: \lieg \to G$ is the 
exponential map and $\epmb = \beta \cdot \epm$.
The space $\xb$ is equipped with two canonical
projection
maps $\proj_1: \xb \to G^{2g}$ and $\proj_2: \xb \to \lieg$, for which
there is the following commutative diagram:

\beq 
\begin{array}{lcr}
\xb  & \stackrel{\proj_2}{\lrar} & \lieg \\
\scriptsize{\proj_1}
\downarrow & \phantom{\stackrel{\proj}{\lrar} } & \downarrow \scriptsize{e_\b} 
 \\
G^{2g}  & \stackrel{\epsr }{\lrar} & G\\ \end{array} \eeq
 A straightforward
argument using the regular value theorem endows 
$\xb$ with a smooth structure on an open dense set $\xb^s$ including 
$\proj_2^{-1}(0)$: see [J1], Proposition 5.4, where an explicit
characterization of the singular locus of $\xb$ is given.

The space $\xb$ carries 
an action of the 
group $G$.
A gauge theoretic construction of 
a $G$-invariant closed 
2-form on $\xb$ 
was given in [J1], and it was shown that
this form is nondegenerate on an open dense set 
in $\xb$ and that the action 
of $G$ is Hamiltonian where the 2-form is nondegenerate.
 There is thus an open dense set in $\xb$ which is a 
smooth finite dimensional Hamiltonian $G$-space such that
the space $\calm_\beta$ is given by symplectic 
reduction of this $G$-space at $0$. 

 By extending the techniques of [W],
 we construct here a 
$G$-invariant closed 
2-form on $\xb$ which is nondegenerate on an open dense set,
 and show that the 
moment map $\mu$ is given by a constant
 multiple of $\proj_2$, as was shown in [J1]. The symplectic form
will in fact turn out to be the same as the one constructed using 
gauge theory (see the Remark at the end of this section): in other
words, one of these symplectic forms is the pullback of the other
under a natural map. It follows from our 
construction that the symplectic quotient  (at 0) 
with respect
to the action of $G$ on $\xb$  is the twisted moduli 
space $\calm_\b$ described above. For $t \in G$, the relative 
character variety $\calm_t$ from Section 4 
is the symplectic reduction of $X$ at the orbit $\calo_\L
\subset \lieg$ (under the adjoint action) that  corresponds
to an element $\L \in \lieg$ for which 
$\exp (\L) = t$. 

To construct the symplectic form, we first construct a form 
$\sigfo \in \Omega^2(\lieg)$ for which 
\beq 
\label{5.1}\epm^* \l = d \sigfo. \eeq 
The existence of such a form
follows from the following standard result (see e.g. 
[Wa], 
Lemma 4.18):
\begin{prop} \label{poincare} [Poincar\'e Lemma]

\noindent (a) If $\g \in \Omega^{p+1} (V)$  where $V$ is a 
vector space, and $d \g = 0 $, then there is a 
form $\sigf   \in \Omega^p(V)$ with 
$\g = d \sigf$. 

\noindent (b) Denote by $\hop$ the map $\Omega^{p+1}
(V) \to \Omega^p(V)$ sending $\g$ to $\sigf$. 
Then $d \hop + \hop d = {\rm id}.$ 
 \end{prop}
 
For $\beta \in \Omega^*(V)$, $\hop \beta$ is given at $v \in V$ by 
\beq 
\label{p.2} (\hop \beta)_v = \int_0^1 F_t^* (\iota_{\bar{v} } \beta) dt \eeq
 where $\bar{v} $ is the vector field on $V$ which takes the constant
value $v$, and $F_t$ is the map $V \to V$ given by multiplication 
by $t$. In our case the form
$$\sigma = \hop (\epm^* \l) $$
is
$G$ -invariant because $\l$ is $G$-invariant and $\epm$ is
a $G$-equivariant map.

We now restrict to the fibre product 
$\xb \subset G^{2g} \times \lieg$.
 For $(h, \L) \in \xb$
we have $\epsr(h) = \epmb(\L)$, so if 
$(H, \z) \in T_h G^{2g} \times \lieg$ represents
an element in the tangent space to $\xb$, we have
\beq \label{n5.5} \epsr_* H = \epmb_* \z. \eeq
We define a 2-form on $\xb$ by 
\beq \tw = \proj_1^* \omega + \proj_2^* \sigfo, \eeq
where $\omega$ was defined in (\ref{4.1}).
We find that 
$$d (\proj_1^* \omega + \proj_2^* \sigfo) = \proj_1^* d \omega + 
\proj_2^* d \sigfo $$
so 
\beq  d \tw (H, \z)  = - \epsr^* \l (H) +  d \sigfo (\z)  
\phantom{aaaa} \mbox{(by (\ref{dsymf}))}
\eeq
$$ \bla = - \l(\epsr_* H) + \l (\epm_* \z) = 
- \l(\epsr_* H) + \l (\epmb_* \z) = 0 . $$
(Here, we have used the fact that $\l$ is invariant under 
multiplication by $\b$, so $\epm^* \l = \epmb^* \l$.)
So we have 

\begin{prop} The 2-form $\tw$ on $G^{2g} \times \lieg$ 
restricts on $\xb$ to a closed form.
\end{prop}

The $G$-invariance of $\tw$ follows because $\sigfo$ and $\omega$
are $G$-invariant and the projection maps $\proj_1$ and $\proj_2$ 
are $G$-equivariant maps.

We may now identify the moment map for the action of 
$G$ on $\xb$: in other words, we find a function 
$\mu: X \to \lieg$ such that 
$\iota_\teta \tw = \eta \cdot d \mu$
where $\teta$ is the vector field on $\xb$ generated by $\eta$
$ \in \lieg$.\footnote{Ordinarily the moment map is specified 
as a map into $\lieg^*$. Here we have used the inner product
to identify $\lieg^*$ with $\lieg$.}
First we recall from (\ref{4.66}) that 
$ (\iota_\teta \omega) = \epsr^* \theta_\eta $.
Now since 
$\sigfo=  \hop (\epm^* \l) $ $= \hop (\epmb^* \l)$, we have 
$$ \iota_\teta \sigfo = \iota_\teta (\hop \epmb^* \l)  = - 
\hop  (\iota_\teta \epmb^* \l)  = - \hop (\epmb^* \iota_\teta \l) $$
\beq \label{5.8}
 \bla  = - \hop (d \epmb^* \theta_\eta)  \onebl  \mbox{(by (\ref{pe4})) }.\eeq
Combining (\ref{5.8}) with 
 Proposition \ref{poincare} (b) we find 
\beq \label{5B} \iota_\teta \sigfo = -\epmb^* \theta_\eta + 
d (I \epmb^* \theta_\eta) \eeq
Adding (\ref{4.66}) and (\ref{5B}) we have 
$$ \iota_\teta (\tw)=  d (\hop \epmb^* \theta_\eta) $$
so that a moment map $\mu: \xb \to \lieg$  for the action of $G$ on 
$\xb$ is  given by 
\beq \eta \cdot \mu = \hop \epmb^* \theta_\eta. \eeq
Now
\beq (\hop \epmb^* \theta_\eta)_\L 
= \int_0^1 F_t^*  \Bigl (\epmb^* \theta_\eta (\L) \Bigr ) 
\onebl \mbox{by (\ref{p.2})} \eeq
$$ \bla = \int_0^1 (\epmb^* \theta_\eta)_{\L t} (\L) = 
\int_0^1 (\theta_\eta)_{\epmb(\L t)} (\epmb_* \L) 
$$ 
$$ = \int_0^1 \eta \cdot (\a + \bar{\a})_{\epmb(\L t)}  (\epmb_* \L) $$ 
$$ = 2 \eta \cdot \L. $$
Thus we have explicitly specified  a moment map for 
the action of $G$, which is equivariant with respect to the adjoint action 
of $G$:
 \begin{prop} \label{p5.3}
A moment map for the 
action of $G$ on $\xb$ is given by the  map $\mu = 
 2 \proj_2: \; (h, \L) \mapsto   2 \L. $ 
\end{prop}

To complete the identification of the form $\tw $   as a 
symplectic form  on an open dense set in 
 $\xb$, one needs the following:
\begin{prop} \label{p5.5}
The form
 $\tw$ is a nondegenerate  bilinear form 
on the Zariski tangent space $T_{(h, \L)} \xb$, for any 
$(h, \L) \in \xb$ for which
$(d \epsr)_h$ is surjective.
\end{prop}
Our proof of this Proposition
 parallels the gauge theory
argument  given in [J1] (see Proposition 3.1 of 
[J1]  for the case $G = SU(2)$). 
This material is 
treated in  [H]  (Theorem 4.4 and Section 5) and [K] (Theorem 4): the proof 
we sketch is essentially
the one given by Huebschmann [H], to whom the group cohomology proof of the
nondegeneracy of the symplectic form on an open neighbourhood of 
the zero locus of the moment map in the 
extended moduli space is due.\footnote{The proof given in [H]
applies to a suitable neighbourhood 
of the zero locus of the moment map in  $\xb$ which is contained in 
$\proj_2^{-1} (\calo_{\rm reg})$, where
$\calo_{\rm reg}$ is the subset of $\lieg$ where the exponential 
map is regular: the subset 
$\proj_2^{-1} (\calo_{\rm reg})$
is  a proper subset of  the smooth locus of 
$\xb$.
We have adapted the 
proof so it applies to the Zariski tangent space $ T_{(h, \L)} \xb$ 
for all points $ (h, \L) \in \xb$ for which
$(d \epsr)_h$ is surjective.}

\nc{\hl}{ {(h, \L) } } 
\noindent{\em Proof of Proposition \ref{p5.5}:}
To establish nondegeneracy of $\tw$
we proceed as follows. Proposition \ref{p5.3} establishes that 
if $\teta$  is the vector field associated to the action 
of $\eta$ on $\xb$, and if $(H, \zeta) \in T_{(h, \L)} \xb$ 
for $(h, \L) \in \xb$ (in other words, 
$(d \epsr)_h H = (d e_\b)_\L \zeta),$ then 
\beq \tw_{(h, \L)}  (\teta, (H,\zeta) ) =  2 \eta \cdot \zeta. \eeq
Thus to establish nondegeneracy of $\tw$ 
at those $\hl$ for which $(d \epsr)_h$ is surjective
(see (\ref{n5.5}))
it suffices to establish it
on the orthocomplement 
in ${\rm Ker} (d \proj_2) \subset T_\hl \xb$ 
of the image of the action of $\lieg$.
This means we must establish that $\omega$ is nondegenerate restricted
to 
\beq \frac{ T_\hl 
\Bigl (\proj_2^{-1} (\L) \subset \xb \Bigr )}{ \{ \teta: 
\eta \in {\rm Stab} (\L) \} }. \eeq

\nc{\gh}{\lieg_{ h} } 
\nc{\geh}{\lieg_{   \epsr(h)}  } 
\nc{\td}{{\delta} }

We have commutative diagrams 
\beq \label{5.4}
\begin{array}{ccccccccc}
\phantom{ 0}  & \phantom{\lrar} & \onebl \onebl 0 \onebl 
 &  \phantom{\lrar} & \phantom{B^1(\FF;} \onebl 0 \phantom{ \gh)} & 
\phantom{ \lrar}
&  \phantom{ B^1 (\ZZ;} \onebl 0 \phantom{ \geh)} & 
\phantom{\lrar} & \phantom{\dots}  \\
\phantom{ 0}  & \phantom{\lrar} & \onebl \downarrow \onebl 
 &  \phantom{\lrar} & \phantom{B^1(\FF;} \downarrow \phantom{ \gh)} & 
\phantom{ \lrar}
&  \phantom{ B^1 (\ZZ;} \downarrow \phantom{ \geh)} & 
\phantom{\lrar} & \phantom{\dots}  \\
 0  & \lrar & K_B &  \lrar & B^1(\FF; \gh) &\stackrel{\td}{  \lrar}
&  B^1 (\ZZ; \geh) & \lrar & \dots  \\
\phantom{ 0}  & \phantom{\lrar} & \onebl \downarrow \onebl 
 &  \phantom{\lrar} & \phantom{B^1(\FF;} \downarrow \phantom{ \gh)} & 
\phantom{ \lrar}
&  \phantom{ B^1 (\ZZ;} \downarrow \phantom{ \geh)} & 
\phantom{\lrar} & \phantom{\dots}  \\
 0  & \lrar & K_C &  \lrar & C^1(\FF; \gh) &  \stackrel{\td}{\lrar}
&  C^1 (\ZZ; \geh) & \lrar & \dots  \\
\phantom{ 0}  & \phantom{\lrar} & \onebl \downarrow \onebl 
 &  \phantom{\lrar} & \phantom{B^1(\FF;} \downarrow \phantom{ \gh)} & 
\phantom{ \lrar}
&  \phantom{ B^1 (\ZZ;} \downarrow \phantom{ \geh)} & 
\phantom{\lrar} & \phantom{\dots}  \\
 0  & \lrar & K_H&  \lrar & H^1(\FF; \gh) &  \stackrel{\td}{\lrar}
&  H^1 (\ZZ; \geh) & \lrar & \dots  \\
\phantom{ 0}  & \phantom{\lrar} & \onebl \downarrow \onebl 
 &  \phantom{\lrar} & \phantom{B^1(\FF;} \downarrow \phantom{ \gh)} & 
\phantom{ \lrar}
&  \phantom{ B^1 (\ZZ;} \downarrow \phantom{ \geh)} & 
\phantom{\lrar} & \phantom{\dots}  \\
\phantom{ 0}  & \phantom{\lrar} & \onebl \onebl 0 \onebl 
 &  \phantom{\lrar} & \phantom{B^1(\FF;} \onebl \onebl 0 \phantom{ \gh)} & 
\phantom{ \lrar}
&  \phantom{ B^1 (\ZZ;} \onebl 0 \phantom{ \geh)} & 
\phantom{\lrar} & \phantom{\dots}  \\
\end{array} \eeq
in which the columns are short exact sequences.
Here, if $\G$ is a discrete 
group (where $\G = \ZZ$ or 
$\G = \FF$) equipped with a 
representation $\rho: \G \to G$, the notation $C^*(\G; \lieg_\rho)$ 
refers to  the Eilenberg-Mac Lane group cochain complex
with coefficients in the $\G$-module $\lieg_\rho$ specified by 
$\rho$ under the adjoint action of $G$ on $\lieg$. 
For $\G = \FF$ or $\G = \ZZ$, we have 
$C^1(\G; \lieg_\rho) = Z^1(\G; \lieg_\rho). $ 
The maps $\td$ are induced by 
$d \epsr$, while $K_B, K_C$ and $ K_H$ are the kernels of the maps $\td$ 
on $B^1, C^1$ and $ H^1$. 
Then the vector space $T_\hl (\proj_2^{-1} (\L) ) $ is 
identified with $K_C \subset C^1 (\FF; \gh) $, while 
$ { \{ \teta: 
\eta \in {\rm Stab} (\L) \} }$ is identified with 
$K_B = K_C \cap B^1 (\FF; \gh)$. Then we see from (\ref{5.4})
 that the quotient
$K_C/K_B$  is canonically identified with $K_H = 
{\rm Ker} \Bigl  (\td: H^1(\FF; \gh) \to H^1 (\ZZ; \geh) \, 
\Bigr ). $

We have the long exact sequence
\beq \label{5.15} 
0 \to H^0 (\FF; \gh) \to H^0 (\ZZ; \geh) \stackrel{\delta^*}{\to}
H^1 (\FF,\ZZ ; \gh) \to H^1 (\FF ; \gh) \to  \eeq
$$ \stackrel{\delta}{\to}
H^1 (\ZZ ; \geh) \to H^2 (\FF,\ZZ ; \gh) \to 0 .  $$
The vector spaces and maps in this sequence satisfy Poincar\'e duality.
Furthermore, the pairing $\, \cdot \,: 
\lieg \otimes \lieg \to \RR$ gives rise to a cup product pairing
\beq 
\label{5.07}
H^1(\FF, \ZZ ; \gh) \otimes H^1(\FF; \gh) \to H^2(\FF, \ZZ; \RR) 
\cong \RR,\eeq and 
Poincar\'e duality in (\ref{5.15}) 
implies that the restriction of this pairing 
to 
$\bigl ( H^1(\FF, \ZZ; \gh)/{\rm Im }(\delta^*) 
\, \bigr ) \otimes \bigl ({\rm Ker}(\delta) 
\subset H^1(\FF; \gh) \bigr  ) $ is nondegenerate. 

Now one may show (see for instance [K] Theorem 4 or 
[H] Theorem 4.4) that this cup product is the restriction 
of the form $\omega$ to ${\rm Ker} (\delta) 
\subset  H^1(\FF; \gh) $ $\cong T_h
\Bigl  (\epsr^{-1} \bigl (C_{\epsr(h)} \bigr  )/G  \Bigr ) $. (Here, 
$C_{\epsr(h)} $ denotes the conjugacy class of $\epsr(h)$ in 
$G$.) 
The nondegeneracy of the  pairing 
arising from the cup product thus completes  the proof of
nondegeneracy.
$\square$ 

\noindent{\em Remark:} Goldman ([G], proof of Proposition 
3.7) has shown that 
$${\rm Im} (d \epsr)_h = \Bigl ( {\rm Lie } (Stab (h) ) \Bigr )^\perp. $$

Suppose $G = SU(n)$ and $\beta$ is a generator of 
$Z(G)$. 
The subset  of $\xb$ where
$(d \epsr)_h$ is surjective then contains the zero level
set of the moment map. 
Thus the following
is an immediate consequence of Proposition \ref{p5.5}:

\begin{corollary} \label{c5.5} Let $G = SU(n)$ and 
suppose $\beta$ is a generator of $Z(G)$. 
Then  $0$ is a regular value of the moment
map
$\mu =  2 \proj_2$. Further, $\calm_\beta$ is a smooth manifold and
the  form $\bar{\omega}$ on $\calm_\beta$ is nondegenerate.\footnote{When
zero is a regular value of the moment map on $\xb$, standard arguments
establish the smoothness of the reduced space $\calm_\beta$.}
\end{corollary}

 \noindent{\em Remark:} At the end of the last Proof, we alluded to 
the identification  of $\omega$ (on quotients of level sets of 
$\proj_2$) with the bilinear form given by the cup product
(\ref{5.07}). It is easy to see that the vector spaces
$ H^1(\FF; \gh) $ and $ H^1(\FF, \ZZ ; \gh) $ (arising from group
cohomology with coefficients in the $\FF$-module $\lieg_h$ 
specified by $h \in {\rm Hom}(\FF, G) $) 
are the same as the vector spaces $H^1(\Sigma-D^2 ; d_A)$ and 
$H^1(\Sigma- D^2, \partial D^2 ; d_A)$  arising in the 
gauge theory description of 
$\xb$ (cf. Section 2.2 of [J1]). 
Here, 
$A$ is a flat connection 
on the punctured surface $\Sigma - D^2$ whose 
holonomy gives rise to the representation $h$ of the fundamental 
group $\FF$.
This identification comes  from the identification  between gauge
equivalence classes of flat $G$ connections and conjugacy classes of 
representations of the fundamental group into $G$, which 
arises from the map sending a flat 
connection to the representation given by its
holonomy.\footnote{Using the holonomy representation, a 
homeomorphism from the extended 
moduli space
(as defined gauge theoretically in Section 2.1 of [J1]) to 
$\xb$ (as defined above) is given in Section 2.3 of [J1]. This 
map identifies the relevant Zariski tangent spaces and gives 
rise to the identification of the symplectic forms.}
Furthermore, the pairing
$ H^1(\FF; \gh) \otimes H^1(\FF, \ZZ ; \gh)  \to \RR$ 
arising from the group cohomology cup product 
is the same as the  pairing
$(\alpha, \beta)  
\mapsto \int_{\Sigma - D^2}  \alpha \cdot 
 \wedge \beta $ that gives the symplectic 
form in the gauge theory description. (Here, 
$\alpha$ and $\beta$ are $d_A$-closed $\lieg$-valued 
1-forms on $\Sigma - D^2$, 
and $\alpha \cdot  \wedge \beta$ is 
the element of $\Omega^2(\Sigma)$ that 
arises from the wedge product 
combined with the pairing 
$ \lieg \otimes \lieg \to \RR$ given by the invariant inner 
product $\, \cdot \,$.)
 Hence the symplectic form we have
constructed on $\xb$ is in  fact the same as the one constructed in [J2] using
gauge theory.

\vspace{0.2in}

{\Large \bf References}

[AB] Atiyah, M.F. and Bott, R., The Yang-Mills equations over Riemann surfaces,
{\em Phil. Trans. R. Soc. Lond.} {\bf A 308} (1982), 523-615.

[B] 
Bott, R., On the Chern-Weil homomorphism and the continuous cohomology
of Lie groups, {\em Advances in Math.} {\bf 11} (1973), 289-303.

[G]  Goldman, W., The symplectic nature of fundamental groups of surfaces,
{\em Advances in Math.} {\bf 54} (1984), 200-225.

[H] Huebschmann, J., Symplectic  and Poisson structures of
certain moduli spaces I, preprint hep-th/9312112 (1993); 
{\em Duke Math. J.} {\bf 80} (1995) 737-756.

[J1] Jeffrey, L.C., Extended moduli spaces of flat connections
on Riemann surfaces, {\em Math. Annalen} {\bf 298} (1994),
667-692.

[J2] Jeffrey, L.C., Group cohomology construction of the cohomology
of moduli spaces of flat connections on 2-manifolds, {\em Duke Math.
J.} {\bf 77} (1995) 407-429.

[JK1] Jeffrey, L.C., Kirwan, F.C., Localization for nonabelian group 
actions, preprint alg-geom/9307001; {\em Topology} 
{\bf 34} (1995) 291-327. 

[JK2] Jeffrey, L.C., Kirwan, F.C., Intersection theory
on moduli spaces of holomorphic bundles on a Riemann 
surface, {\em Elec. Res. Notices AMS} 
{\bf 1} (No. 2) (1995) 57-71.

[JW1] Jeffrey, L.C., Weitsman, J. Toric structures on the 
moduli space of flat connections on a Riemann surface: volumes 
and the moment map, {\em Advances in Math.} {\bf 109} (1994),
151-168.

[JW2] Jeffrey, L.C., Weitsman, J. Torus actions and the topology and
symplectic geometry of flat connections on 2-manifolds,
 in  {\em Topology, Geometry and 
Field Theory} (Proceedings of the Taniguchi Foundation Symposium on 
Low Dimensional Topology and Topological Field Theory, Kyoto, 
January 1993), ed. K. Fukaya, M. Furuta, T. Kohno, 
D. Kotschick, World Scientific,
1994.

[K] 
Karshon, Y., An algebraic proof for the symplectic structure of moduli space,
{\em Proc. Amer. Math. Soc.} {\bf 116} (1992), 591-605.

[KT] Karshon, Y., Tolman, S., The moment map and line bundles over 
presymplectic toric manifolds, {\em J. Diff. Geom.}
{\bf 38} (1993), 465-484.

[MS] Mehta, V., Seshadri, C.S., Moduli of vector bundles on curves
with parabolic structure, {\em Math. Annalen} {\bf 248} (1980), 205-239.

[Se] C.S. Seshadri, {\em Fibr\'es Vectoriels sur les 
Courbes Alg\'ebriques,} {\em Ast\'erisque} {\bf 96} (1982).

[Sh]
Shulman, H.B., Characteristic classes and foliations, Ph.D. Thesis,
University of California, Berkeley (1972).

[Wa] Warner, F.W., {\em Foundations of Differentiable Manifolds
and Lie Groups}, Springer-Verlag, 1983.

[W] Weinstein, A., The symplectic structure on moduli space, 
in {\em The  Floer Memorial Volume}, ed. H. Hofer, C. Taubes, 
A. Weinstein,
E. Zehnder, Birkh\"auser (Progress in Math. {\bf 133}) (1995) 
627-635.

[Wi]  E. Witten, { Two dimensional gauge theories
revisited}, preprint hep-th/9204083;
 {\em J. Geom. Phys.} {\bf 9} (1992), 303-368.

\end{document}